\documentclass[11pt]{article}

\usepackage{amsmath,makeidx,amssymb,amsthm,amscd,a4}
\usepackage{graphicx,epsfig,latexsym,bbm,index}
\usepackage[active]{srcltx}

\newtheorem{theor}{Theorem}
\newtheorem{lemma}[theor]{Lemma}

\begin{document}

\title{Semiconductors and Dirichlet-to-Neumann maps}

\author{Antonio Leit\~ao \\
\small Department of Mathematics, Federal University of St.\,Catarina, \\
\small 88040-900 Florianopolis, Brazil ({\tt aleitao@mtm.ufsc.br}) }

\date{}

\maketitle

\begin{abstract}
We investigate the problem of identifying discontinuous doping profiles in
semiconductor devices from data obtained by the stationary voltage-current
(VC) map. The related inverse problem correspond to the inverse problem for
the Dirichlet-to-Neumann (DN) map with partial data.
\end{abstract}

\section{Introduction} \label{sec:0}

In this paper we investigate the problem of identifying discontinuous doping
profiles in semiconductor devices from data obtained by the VC map for the
linearized stationary bipolar model (close to equilibrium).
Two different methods of data acquisition are considered:

1) Current flow measurements through a contact;

2) Pointwise measurements of the current density.

\noindent
The related inverse problems correspond to the inverse problem for the DN
map with partial data.

We propose a framework to handle the inverse problems and analyze relevant
properties of the parameter-to-output maps. Moreover, we present a numerical
experiment for the case of pointwise measurements of the current density.

The paper is outlined as follows: In Section~\ref{sec:1} we present the
transient end stationary drift diffusion equations. In the Section~\ref{sec:2}
we introduce the VC map and derive the underlying model for the analysis
presented in this paper, namely the {\em linearized stationary bipolar case
(close to equilibrium)}. Two inverse problems corresponding to different data
acquisition procedures are introduced in Section~\ref{sec:3}.
Some regularity properties of the related parameter-to-output maps are
verified in this section. In Section~\ref{sec:4} we present some numerical
results for a level set type iterative method and pointwise measurements of
the current density. This experiment indicate that a single measurement may
suffices to identify the doping profile.

\section{Drift diffusion equations} \label{sec:1}

\paragraph{The transient model}
\mbox{} \medskip

\noindent
The basic semiconductor device equations consist of the Poisson equation
(\ref{eq:VanR1}), the continuity equations for electrons (\ref{eq:VanR2})
and holes (\ref{eq:VanR3}), and the current relations for electrons
(\ref{eq:VanR4}) and holes (\ref{eq:VanR5}).
\begin{subequations} \label{eq:VanR} \begin{eqnarray}
{\rm div} (\epsilon \nabla V) &       \hskip-1.2cm \label{eq:VanR1}
   =  q(n - p - C) & {\rm in}\ \Omega \times (0,T) \\
{\rm div}\, J_n &                     \hskip-1.5cm \label{eq:VanR2}
   =  q ( \partial_t n + R) & {\rm in}\ \Omega \times (0,T) \\
{\rm div}\, J_p &                     \hskip-1.25cm \label{eq:VanR3}
   =  q (-\partial_t p - R) & {\rm in}\ \Omega \times (0,T)  \\
J_n &                                 \hskip-0.15cm \label{eq:VanR4}
   =  q ( D_n \nabla n - \mu_n n \nabla V) & {\rm in}\ \Omega \times (0,T) \\
J_p &                                 \label{eq:VanR5}
   =  q (-D_p \nabla p - \mu_p p \nabla V) & {\rm in}\ \Omega \times (0,T) \,.
\end{eqnarray} \end{subequations}
This system is defined in $\Omega \times (0,T)$, where $\Omega \subset
\mathbb{R}^d$ ($d=1,2,3$) is a domain representing the semiconductor
device.
Here $V$ denotes the electrostatic potential ($- \nabla V$ is the electric
field, $E=\vert\nabla V\vert$)), $n$ and $p$ are the concentration of free
carriers of negative charge (electrons) and positive charge (holes)
respectively and $J_n$ and $J_p$ are the densities of the electron and the
hole current respectively. $D_n$ and $D_p$ are the diffusion coefficients for
electrons and holes respectively. $\mu_n$ and $\mu_p$ denote the mobilities of
electrons and holes respectively. The positive constants $\epsilon$ and $q$
denote the permittivity coefficient (for silicon) and the elementary charge
(see Appendix).

The function $R$ has the form $R = \mathcal{R}(n,p,x) (np - n_i^2)$ and
denotes the {\em recombination-generation rate} ($n_i$ is the intrinsic
carrier density). The {\em bandgap} is relatively large for semiconductors
(gap between valence and conduction band), and a significant amount of energy
is necessary to transfer electrons from the valence and to the conduction
band. This process is called generation of electron-hole pairs. On the
other hand, the reverse process corresponds to the transfer of a
conduction electron into the lower energetic valence band. This process
is called recombination of electron-hole pairs. In our model these
phenomena are described by the recombination-generation rate $R$.
Frequently adopted in the literature are the Shockley Read Hall model
($\mathcal{R}_{SRH}$) and the Auger model ($\mathcal{R}_{AU}$), defined by
$$  \mathcal{R}_{SRH} := [ \tau_p(n+n_i) + \tau_p(p+n_i) ]^{-1}
    \, , \ \ \mathcal{R}_{AU} := C_n n + C_p p \, , $$
where $C_n$, $C_p$, $\tau_n$ and $\tau_p$ are positive constants
(see Appendix).

The function $C(x)$ models a preconcentration of ions in the crystal, so $C(x)
= C_{+}(x) - C_{-}(x)$ holds, where $C_{+}$ and $C_{-}$ are concentrations of
negative and positive ions respectively. In those subregions of $\Omega$ in
which the preconcentration of negative ions predominate (P-regions), we have
$C(x) < 0$. Analogously, we define the N-regions, where $C(x) > 0$ holds.
The boundaries between the P-regions and N-regions (where $C$ change sign)
are called {\em pn-junctions}.

In the sequel we turn our attention to the boundary conditions.
We assume the boundary $\partial\Omega$ of $\Omega$ to be divided into two
nonempty disjoint parts: $\partial\Omega = \overline{\partial\Omega_N} \cup
\overline{\partial\Omega_D}$. The Dirichlet part of the boundary
$\partial\Omega_D$ models the Ohmic contacts, where the potential $V$ as
well as the concentrations $n$ and $p$ are prescribed. The Neumann part
$\partial\Omega_N$ of the boundary corresponds to insulating surfaces,
thus a zero current flow and a zero electric field in the normal
direction are prescribed.
The Neumann boundary conditions for system (\ref{eq:VanR1}) -- (\ref{eq:VanR5})
read:
\begin{equation} \label{eq:VanR-bcN}
 \frac{\partial V}{\partial\nu}(x,t) =
   \frac{\partial n}{\partial\nu}(x,t) =
   \frac{\partial p}{\partial\nu}(x,t) = 0 \, , \
   \partial\Omega_N \times [0,T] \, .
\end{equation}
Moreover, at $\partial\Omega_D \times [0,T]$, the following Dirichlet boundary
conditions are imposed:
\begin{subequations} \label{eq:VanR-bcD} \begin{eqnarray}
V(x,t) &\!\!\! = &\!\!\! V_D(x,t)
       \, = \, U(x,t) + V_{\rm bi}(x)
       \, = \, U(x,t) + U_T \, \ln \left( \frac{n_D(x)}{n_i} \right) \\[1ex]
n(x,t) &\!\!\! = &\!\!\! n_D(x)
       \, = \, \frac{1}{2} \left(C(x) + \sqrt{C(x)^2 + 4 n_i^2}\right) \\[1ex]
p(x,t) &\!\!\! = &\!\!\! p_D(x)
       \, = \, \frac{1}{2} \left(-C(x) + \sqrt{C(x)^2 + 4 n_i^2}\right) \, .
\end{eqnarray} \end{subequations}
Here the function $U(x,t)$ denotes the applied potential, the constant
$U_T$ represents the thermal voltage, and $V_{\rm bi}$ is given logarithmic
function \cite{BELM04}.
We shall consider the simple situation $\partial\Omega_D = \Gamma_0 \cup
\Gamma_1$, which occurs, e.g., in a diode. The disjoint boundary parts
$\Gamma_i$, $i=0,1$, correspond to distinct contacts.
Differences in $U(x)$ between different segments of $\partial\Omega_D$
correspond to the applied bias between these two contacts. Moreover,
the initial conditions $n(x,0) \ge 0$, $p(x,0) \ge 0$ have to be imposed.

\paragraph{The stationary model}
\mbox{} \medskip

\noindent
In this paragraph we turn our attention to the stationary drift diffusion
equations. We disconsider the thermal effects and assume further
$\frac{\partial n}{\partial t} = \frac{\partial n}{\partial t} = 0$.
Thus, the {\em stationary drift diffusion model} is derived from
(\ref{eq:VanR1}) -- (\ref{eq:VanR5}) in a straightforward way.
Next, motivated by the Einstein relations $D_n = U_T \mu_n$ and
$D_p = U_T \mu_p$ (a standard assumption about the mobilities and
diffusion coefficients), one introduces the so-called
{\em Slotboom variables} $u$ and $v$, which are related to the
original $n$ and $p$ variables by the formula:
\begin{equation} \label{eq:slotboom}
n(x) = n_i \exp\left(\frac{ V(x)}{U_T}\right)\, u(x)\, , \ \ \
p(x) = n_i \exp\left(\frac{-V(x)}{U_T}\right)\, v(x) \, .
\end{equation}
For convenience, we rescale the potential and the mobilities, i.e.
$ V(x) \ \leftarrow \ V(x) / U_T$, \,$\mu_n \leftarrow q U_T \mu_n$,\,
$\mu_p \leftarrow q U_T \mu_p$. It is obvious to check that the current
relations now read $J_n = \mu_n n_i \, e^{ V} \nabla u$,\,
$J_p = -\mu_p n_i \, e^{-V} \nabla v$.

Next we write the stationary drift diffusion equations in terms of $(V,u,v)$
\begin{subequations} \label{eq:dd-sys-nlin} \begin{eqnarray}
\lambda^2 \, \Delta V & \hskip-0.15cm \label{eq:dd-sys1}
  = \ \delta^2 \big(e^Vu - e^{-V}v\big) - C(x) & {\rm in}\ \Omega \\
{\rm div}\, J_n & \hskip-0.3cm \label{eq:dd-sys2}
  = \ \delta^4 \, Q(V,u,v,x) \, (u v - 1)      & {\rm in}\ \Omega \\
{\rm div}\, J_p & \label{eq:dd-sys3}
  = \ - \delta^4 \, Q(V,u,v,x) \, (u v - 1)    & {\rm in}\ \Omega \\
V & \hskip-1.62cm = \ V_D \ = \ U + V_{\rm bi} & {\rm on}\ \partial\Omega_D
\label{eq:dd-sys4} \\
u & \hskip-2.2cm = \ u_D \ = \ e^{-U}          & {\rm on}\ \partial\Omega_D
\label{eq:dd-sys5} \\
v & \hskip-2.45cm = \ v_D \ = \ e^{U}           & {\rm on}\ \partial\Omega_D
\label{eq:dd-sys6} \\
\nabla V \cdot \nu & \hskip-0.6cm = \ J_n\cdot\nu \ = \ J_p\cdot\nu \ = \ 0
                   & \rm on \ \partial\Omega_N \, , \label{eq:dd-sys7}
\end{eqnarray} \end{subequations}
where $\lambda^2 := \epsilon/(q U_T)$ is the Debye length of the device,
$\delta^2 := n_i$ and the function $Q$ is defined implicitly by the
relation $Q(V,u,v,x) = {\cal R}(n,p,x)$.%
\footnote{Notice the applied potential has also to be rescaled:
$U(x) \leftarrow U(x) / U_T$.}

One should notice that, due to the thermal equilibrium assumption, it follows
$np = n_i^2$, and the assumption of vanishing space charge density gives
$n-p-C = 0$, for $x \in \partial\Omega_D$. This fact motivates the
boundary conditions on the Dirichlet part of the boundary.

It is worth mentioning that, in a realistic model, the mobilities
$\mu_n$ and $\mu_p$ usually depend on the electric field strength
$|\nabla V|$.
In what follows, we assume that $\mu_n$ and $\mu_p$ are positive constants.
This assumption simplifies the subsequent analysis, allowing us to concentrate
on the inverse doping problems. As a matter of fact, this dependence could
be incorporated in the model without changing the results described in the
sequel.

Existence and uniqueness of solutions for system (\ref{eq:dd-sys-nlin})
can only be guaranteed for small applied voltages. Therefore, it is
reasonable to consider, instead of this system, its linearized version
around the equilibrium point $U \equiv 0$. We shall return to this point
in the next section, where the VC map is introduced.

\section{A simplified model} \label{sec:2}

In the sequel we make some simplifying assumptions on the stationary drift
diffusion equations introduced in Section~\ref{sec:1} and derive a special
case which will serve as underlying model for the inverse problem investigated
in Section~\ref{sec:3}.

\paragraph{The linearized stationary drift diffusion equations
(close to equilibrium)}
\mbox{} \medskip

\noindent
We begin this paragraph by introducing the {\em thermal equilibrium}
assumption for the stationary drift diffusion equations. This is a previous
step to derive a linearized system of stationary drift diffusion equations
(close to equilibrium).

The thermal equilibrium assumption refers to the condition in which the
semiconductor is not subject to external excitations, except for a uniform
temperature, i.e. no voltages or electric fields are applied.
It is worth noticing that, under the thermal equilibrium assumption, all
externally applied potentials to the semiconductor contacts are zero
(i.e. $U(x) = 0$). Moreover, the thermal generation is perfectly balanced
by recombination (i.e. $\mathcal R = 0$).

If the applied voltage satisfies $U = 0$, one immediately sees that the
solution of system (\ref{eq:dd-sys1}) -- (\ref{eq:dd-sys7}) simplifies to
$(V,u,v) = (V^0, 1,1)$, where $V^0$ solves
\begin{subequations}  \label{eq:equil-case} \begin{eqnarray}
\lambda^2 \, \Delta V^0 &                       \label{eq:equil-caseA}
   = \ e^{V^0} - e^{-V^0} - C(x)              & {\rm in}\ \Omega \\
V^0 & \hskip-2.1cm                              \label{eq:equil-caseB}
   = \ V_{\rm bi}(x)                          & {\rm on}\ \partial\Omega_D \\
\nabla V^0 \cdot \nu & \hskip-2.95cm            \label{eq:equil-caseC}
   = \ 0                                      & {\rm on}\ \partial\Omega_N \, .
\end{eqnarray} \end{subequations}

In the bipolar model discussed below we shall be interested in the linearized
drift diffusion system at the equilibrium. Keeping this in mind, we compute
the Gateaux derivative of the solution of system (\ref{eq:dd-sys1}) --
(\ref{eq:dd-sys7}) with respect to the voltage $U$ at the point $U \equiv 0$
in the direction $h$. This directional derivative is given by the solution
$(\hat V, \hat u, \hat v)$ of
\begin{subequations}  \label{eq:dd-sys-lin} \begin{eqnarray}
\lambda^2 \, \Delta \hat V & 
  = \ e^{V^0} \hat u + e^{-V^0} \hat v + ( e^{V^0} + e^{-V^0} ) \hat V
                                               & {\rm in}\ \Omega \\
{\rm div}\, (\mu_n e^{V^0} \nabla \hat u) & \hskip-2.38cm
  = \ Q_0(V^0,x) (\hat u + \hat v)             & {\rm in}\ \Omega \\
{\rm div}\, (\mu_p e^{-V^0} \nabla \hat v) & \hskip-2.38cm
  = \ Q_0(V^0,x) (\hat u + \hat v)             & {\rm in}\ \Omega \\
\hat V & \hskip-5.0cm
  = \ h                                        & {\rm on}\ \partial\Omega_D \\
\hat u & \hskip-4.7cm
  = \ -h                                       & {\rm on}\ \partial\Omega_D \\
\hat v & \hskip-5.0cm
  = \  h                                       & {\rm on}\ \partial\Omega_D \\
\nabla V^0 \cdot\nu & \hskip-1.7cm
  = \ \nabla\hat u \cdot\nu \, = \, \nabla\hat v \cdot\nu \, = \, 0
                                              & {\rm on}\ \partial\Omega_N \, ,
\end{eqnarray} \end{subequations}
where the function $Q_0$ satisfies $Q_0(V^0,x) = Q(V^0,1,1,x)$.

\paragraph{Linearized stationary bipolar case (close to equilibrium)}
\mbox{} \medskip

\noindent
In this paragraph we present a special case, which plays a key rule in the
modeling inverse doping problems related to {\em current flow} measurements.

The following discussion is motivated by the stationary VC map
$$ \begin{array}{rcl}
   \Sigma_C: H^{3/2}(\partial\Omega_D) & \to & \mathbb R \, . \\
   U & \mapsto & \displaystyle\int_{\Gamma_1} (J_n+J_p)\cdot\nu \, ds
   \end{array} $$
Here $(V,u,v)$ is the solution of system (\ref{eq:dd-sys-nlin}) for an
applied voltage $U$. This operator models practical experiments where
{\em voltage-current data} are available, i.e. measurements of the
averaged outflow current density on $\Gamma_1 \subset \partial\Omega_D$.

The {\em linearized stationary bipolar case (close to equilibrium)}
corresponds to the model obtained from the drift diffusion equations
(\ref{eq:dd-sys-nlin}) by linearizing the VC map at $U \equiv 0$. This
simplification is motivated by the fact that, due to hysteresis effects for
large applied voltage, the VC map can only be defined as a single-valued
function in a neighborhood of $U=0$.
Moreover, the following simplifying assumptions are also taken into account:

\begin{itemize}
\item[{\it A1)}] \ The electron mobility $\mu_n$ and hole mobility $\mu_p$
are constant;
\item[{\it A2)}] \ No recombination-generation rate is present, i.e.
${\cal R} = 0$ (or $Q_0 = 0$).
\end{itemize}

An immediate consequence of our assumptions is the fact that the 
Poisson equation and the continuity equations decouple. Indeed,
from (\ref{eq:dd-sys-lin}) we see that the Gateaux derivative of
the VC map $\Sigma_C$ at the point $U=0$ in the direction
$h \in H^{3/2}(\partial\Omega_D)$ is given by the expression
\begin{equation} \label{eq:def-sigma-prime-C}
\Sigma'_C(0) h = \int_{\Gamma_1}
\left( \mu_n \, e^{V_{\rm bi}} \hat{u}_\nu -
\mu_p \, e^{-V_{\rm bi}} \hat{v}_\nu \right) \, ds ,
\end{equation}
where $(\hat{u}, \hat{v})$ solve
\begin{subequations}  \label{eq:bipol-stat}  \begin{eqnarray}
{\rm div}\, (\mu_n e^{V^0} \nabla \hat{u})   & \hskip-1.6cm
 = \ 0                \label{eq:bipol-statA} & {\rm in}\ \Omega \\
{\rm div}\, (\mu_p e^{-V^0} \nabla \hat{v})  & \hskip-1.6cm
 = \ 0                \label{eq:bipol-statB} & {\rm in}\ \Omega \\
\hat{u} & \hskip-1.3cm
 = \ -h                                      & {\rm on}\ \partial\Omega_D \\
\hat{v} & \hskip-1.55cm
 = \ h                                       & {\rm on}\ \partial\Omega_D \\
\nabla\hat{u} \cdot\nu &
 = \ \nabla\hat{v} \cdot\nu \ = \ 0          & {\rm on}\ \partial\Omega_N
\end{eqnarray} \end{subequations}
and $V^0$ is the solution of the equilibrium problem (\ref{eq:equil-case});
see Lemma~\ref{prop:bemp31}.

Notice that the solution of the Poisson equation can be computed a priori,
since it does not depend on $h$.
The application $\Sigma_C'(0)$ maps the Dirichlet data for $(\hat u,\hat v)$
to a weighted sum of their Neumann data and can be compared with the
DN map in the {\em Electrical Impedance Tomography} (EIT).

\section{Inverse Problems for Semiconductors} \label{sec:3}

We begin this section verifying that the stationary VC map $\Sigma_C$,
introduced in Section~\ref{sec:2}, is well defined in a suitable
neighborhood of $U=0$.

\begin{lemma} \mbox{\bf \cite[Proposition 3.1]{BEMP01}}
\label{prop:bemp31}
In (\ref{eq:dd-sys-nlin}),
for each applied voltage $U \in B_r(0)$ $\subset H^{3/2}(\partial\Omega_D)$
with $r>0$ sufficiently small, the current $J \cdot \nu \in H^{1/2}
(\Gamma_1)$ is uniquely defined. Furthermore, $\Sigma_C:
H^{3/2}(\partial\Omega_D) \to H^{1/2}(\Gamma_1)$ is continuous
and continuously differentiable in $B_r(0)$. Moreover, it's derivative
in direction $h \in H^{3/2}(\partial\Omega_D)$ is given by the operator
$\Sigma'_C(0)$ defined in (\ref{eq:def-sigma-prime-C}).
\end{lemma}

As a matter of fact, we can actually prove that, since $(\hat u,\hat v)$
in (\ref{eq:bipol-stat}) depend continuously (in $H^2(\Omega)^2$) on the
boundary data $U \in H^{3/2}(\partial\Omega_D)$, it follows from the
boundedness and compactness of the trace operator $\gamma: H^2(\Omega)
\to H^{1/2}(\Gamma_1)$ that $\Sigma'_C(0)$ is a bounded and compact
operator.

Lemma~\ref{prop:bemp31} establishes a basic property to consider
the inverse problem of reconstructing the doping profile $C$ from the
VC map. In the sequel we shall consider two possible inverse problems
for this map.

\paragraph{Current flow measurements through a contact}
\mbox{} \medskip

\noindent
In this first inverse problem we assume that, for each $C$, the output
is  given by $\Sigma_C'(0) U_j$ for some $U_j$. A realistic experiment
corresponds to measure, for given $\{ U_j \}_{j=1}^N$, with $\|U_j\|$
small, the outputs
$$ \big\{ \Sigma_C'(0) U_j\ |\ \ j=1,\cdots,N \big\} $$
(recall that $\Sigma_C(0) = (V^0,1,1)$).
In practice, the functions $U_j$ are chosen to be piecewise constant on
the contact $\Gamma_1$ and to vanish on $\Gamma_0$. 
From the definition of $\Sigma_C'(0)$ we derive the following abstract
formulation of the inverse doping profile problem for the VC map:
\begin{equation} \label{eq:ip-abstract}
 F(C) \ = \ Y \, ,
\end{equation}
where
\begin{enumerate}
\item[1)] $\{ U_j \}_{j=1}^N \subset H^{3/2}(\partial\Omega_D)$ are
fixed voltage profiles satisfying $U_j |_{\Gamma_1} = 0$;
\item[2)] Parameter: \ $C = C(x) \ \in \ L^2(\Omega) =: \mathcal X$;
\item[3)] Output: \ $Y = \big\{ \Sigma_C'(0) U_j \big\}_{j=1}^N \in
\mathbb R^N =: \mathcal Y$;
\item[4)] Parameter-to-output map: \ $F: \mathcal X \to \mathcal Y$.
\end{enumerate}
The domain of definition of the operator $F$ is
$$ D(F) := \{ C \in L^\infty(\Omega) ; \, C_m \le C(x) \le C_M,
   \mbox{ a.e. in } \Omega \} \, , $$
where $C_m$ and $C_M$ are suitable positive constants.

The inverse problem described above corresponds to the problem of identifying
the doping profile $C$ from the linearized stationary VC map at $U=0$
(see bipolar case in Section~\ref{sec:2}).

The non-linear parameter-to-output operator $F$ is well defined and
Fr\'echet differentiable in its domain of definition $D(F)$. This
assertion follows from standard regularity results in PDE theory
(see, e.g., \cite[Propositions~2.2 and~2.3]{BELM04}).

It is worth noticing that the solution of the Poisson equation can be
computed {\em a priori}. The remaining problem (coupled system
(\ref{eq:bipol-stat}) for $(\hat u,\hat v)$) is quite similar to the
problem of EIT. In this inverse problem the aim is to identify the
conductivity $q = q(x)$ in the equation:
$$ -{\rm div}\,(q \nabla u) \ = \ f \ \  {\rm in}\ \Omega \, , $$
from measurements of the {\em Dirichlet-to-Neumann map}, which maps
the applied voltage $u|_{\partial\Omega}$ to the electrical flux
$q u_\nu|_{\partial\Omega}$. The application $\Sigma_C'(0)$ maps the
Dirichlet data for $\hat{u}$ and $\hat{v}$ to the weighted sum of their
Neumann data. It can be seen as the counterpart of electrical impedance
tomography for common conducting materials.

\paragraph{Pointwise measurements of the current density}
\mbox{} \medskip

\noindent
In the sequel, we investigate a different inverse problem related to the
VC map. Differently from the previous paragraph, we shall assume that the
VC operator maps the Dirichlet data for $\hat{u}$ and $\hat{v}$ in
(\ref{eq:bipol-stat}) to the sum of their Neumann data, i.e.
$$ \begin{array}{rcl}
     \Sigma_C: H^{3/2}(\partial\Omega_D) & \to & H^{1/2}(\Gamma_1) \\
     U & \mapsto & (J_n + J_p)\cdot\nu |_{\Gamma_1}
   \end{array} $$
where functions $V$, $\hat u$, $\hat v$, $J_n$, $J_p$ and $U$ have the
same meaning as in Section~\ref{sec:2}. It is immediate to observe that
the Gateaux derivative of the VC map $\Sigma_C$ at the point $U=0$ in the
direction $h \in H^{3/2}(\partial\Omega_D)$ is given by
\begin{equation}
\Sigma'_C(0) h = \left( \mu_n \, e^{V_{\rm bi}} \hat{u}_\nu -
\mu_p \, e^{-V_{\rm bi}} \hat{v}_\nu \right) |_{\Gamma_1} \, ,
\end{equation}
where $(\hat{u}, \hat{v})$ solve system (\ref{eq:bipol-stat}). Notice that,
for each applied voltage $U$, the VC map associates a scalar valued function
defined on $\Gamma_1$. In this case, the outputs $\Sigma'_C(0) U_j$ are in
a data space which is larger than in the case of current flow measurements.

Again we can derive an abstract formulation of type (\ref{eq:ip-abstract})
for the inverse doping profile problem for the linearized stationary VC map
with pointwise measurements of the current density. The only difference to the
framework described in the previous paragraph concerns the definition of the
Hilbert space $Y$, which is now defined by:
\begin{enumerate}
\item[3')] Output: \ $Y = \big\{ \Sigma_C'(0) U_j \big\}_{j=1}^N \in
L^2(\Gamma_1)^N =: \mathcal Y$;
\end{enumerate}
The domain of definition of the operator $F$, remains unaltered.

It is immediate to observe that the model concerning current flow measurements
carries less information about the unknown parameter than the model related to
pointwise measurements does.
In so far, the inverse problem related with the first measurement type is
harder to solve.

\section{A Numerical Experiment} \label{sec:4}

In this section we apply numerical methods to solve an inverse doping profile
problem related to the VC map. We consider the stationary linearized bipolar
model with pointwise measurements of the current density.

In the sequel we consider the bipolar model introduced in Section~\ref{sec:2}.
It follows from the assumption $Q = 0$ that the Poisson equation
(\ref{eq:equil-caseA}) and the continuity equations (\ref{eq:bipol-statA}),
(\ref{eq:bipol-statB}) decouple.
The inverse doping profile problem corresponds to the identification of
$C = C(x)$ from pointwise measurements of the total current density
$J$ at the contact $\Gamma_1$, namely
$$ J |_{\Gamma_1} = (J_n + J_p) |_{\Gamma_1} =
   ( \mu_n e^{V_{\rm bi}} \hat{u}_\nu - \mu_p e^{-V_{\rm bi}} \hat{v}_\nu )
|_{\Gamma_1} $$
(compare with the Gateaux derivative of the VC map $\Sigma_C$ at the point
$U=0$ in (\ref{eq:def-sigma-prime-C})).
Here $(V^0, \hat u, \hat v)$ solve, for each applied voltage $U$, the system
(\ref{eq:equil-case}), (\ref{eq:bipol-stat}), with $h$ substituted by $U$.

Notice that we can split the inverse problem in two parts:
First we define the function $\gamma(x) := e^{V^0(x)}$, $x \in \Omega$,
and solve the parameter identification problem
\begin{equation} \label{eq:num-d2nB}
\begin{array}{r@{\ }c@{\ }l@{\ }l}
   {\rm div}\, (\mu_n \gamma \nabla \hat u) & = & 0 &  {\rm in}\ \Omega \\
   \hat u & = & - U(x) & {\rm on}\ \partial\Omega_D \\
   \nabla \hat u \cdot \nu & = & 0 & {\rm on}\ \partial\Omega_N
\end{array}
\hskip0.6cm
\begin{array}{r@{\ }c@{\ }l@{\ }l}
   {\rm div}\, (\mu_p \gamma^{-1} \nabla \hat v) & = & 0 &  {\rm in}\ \Omega \\
   \hat v & = & U(x) & {\rm on}\ \partial\Omega_D \\
   \nabla \hat v \cdot \nu & = & 0 & {\rm on}\ \partial\Omega_N
\end{array}
\end{equation}
for $\gamma$, from measurements of 
$( \mu_n \gamma \hat{u}_\nu - \mu_p \gamma^{-1} \hat{v}_\nu )
|_{\Gamma_1}$. The second step consists in the determination of $C$ in
$$ C(x) \ = \ \gamma(x) - \gamma^{-1}(x) - \lambda^2 \,
   \Delta (\ln \gamma(x)) \, ,\ x \in \Omega \, . $$
Since the evaluation of $C$ from $\gamma$ can be explicitely performed
in a stable way, we shall focus on the problem of identifying the function
parameter $\gamma$ in (\ref{eq:num-d2nB}).

Summarizing, the inverse doping profile problem in the linearized bipolar
model for pointwise measurements of the current density reduces to the
identification of the parameter $\gamma$ in (\ref{eq:num-d2nB}) from
measurements of the DN map
$$  \Lambda_\gamma : \begin{array}[t]{rcl}
    H^{3/2}(\partial\Omega_D) & \to & H^{1/2}(\Gamma_1) \, . \\
    U & \mapsto & ( \mu_n \gamma \hat{u}_\nu - \mu_p \gamma^{-1} \hat{v}_\nu )
                  |_{\Gamma_1}
    \end{array} $$

If we take into account the restrictions imposed by the practical experiments
described in Section~\ref{sec:3}, it follows:
\begin{itemize}
\item[{\em i)}] The voltage profiles $U \in H^{3/2}(\partial\Omega_D)$ must
satisfy $U |_{\Gamma_1} = 0$;
\item[{\em ii)}] The identification of $\gamma$ has to be performed from a
finite number of measurements, i.e. from the data
\begin{equation} \label{eq:data-bp-pm}
\big\{ (U_j, \Lambda_\gamma(U_j)) \big\}_{j=1}^N
   \in \big[ H^{3/2}(\Gamma_0) \times H^{1/2}(\Gamma_1) \big]^N .
\end{equation}
\end{itemize}

\begin{figure}[ht]
\centerline{
\includegraphics[width=6.6cm]{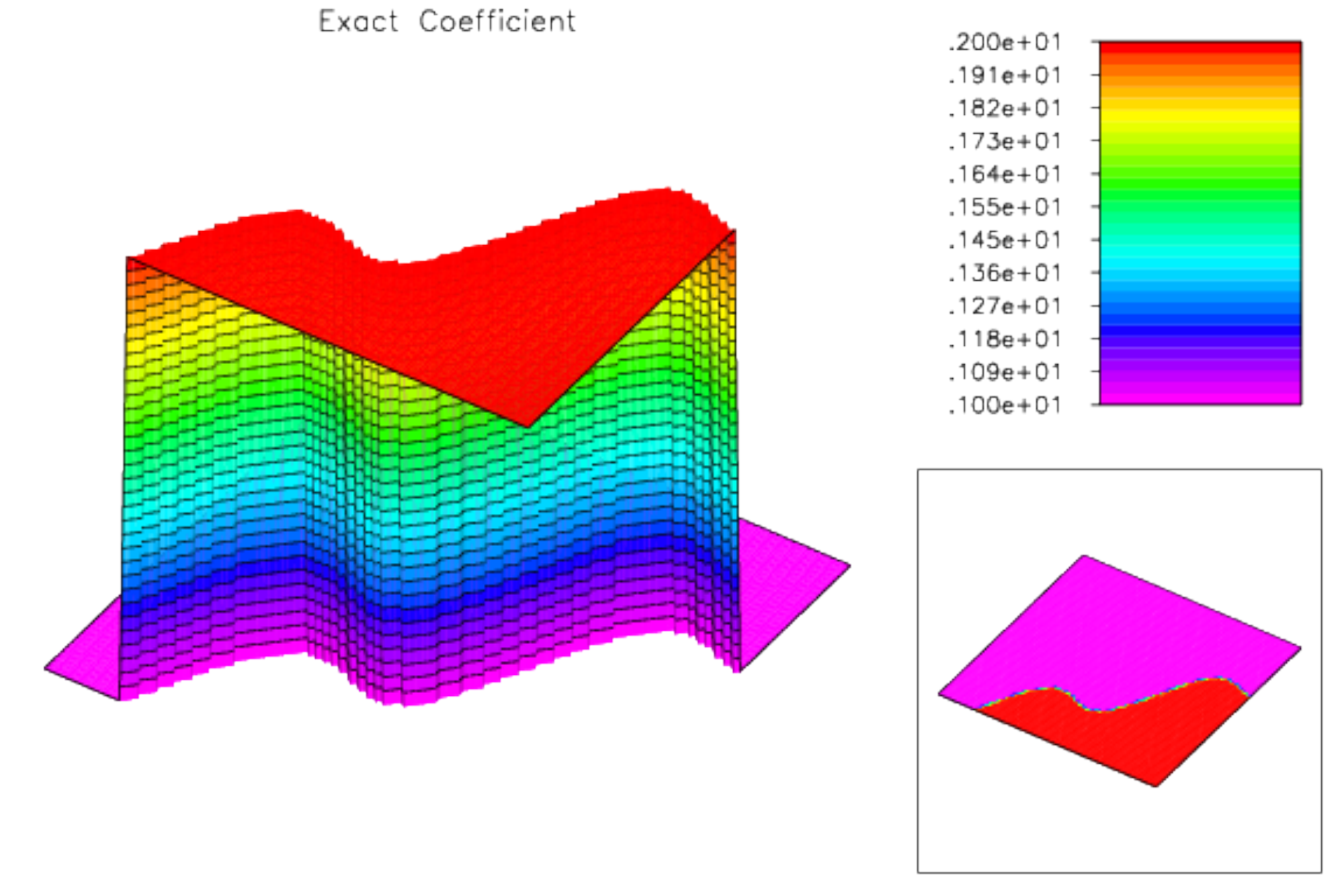}
\includegraphics[width=6.6cm]{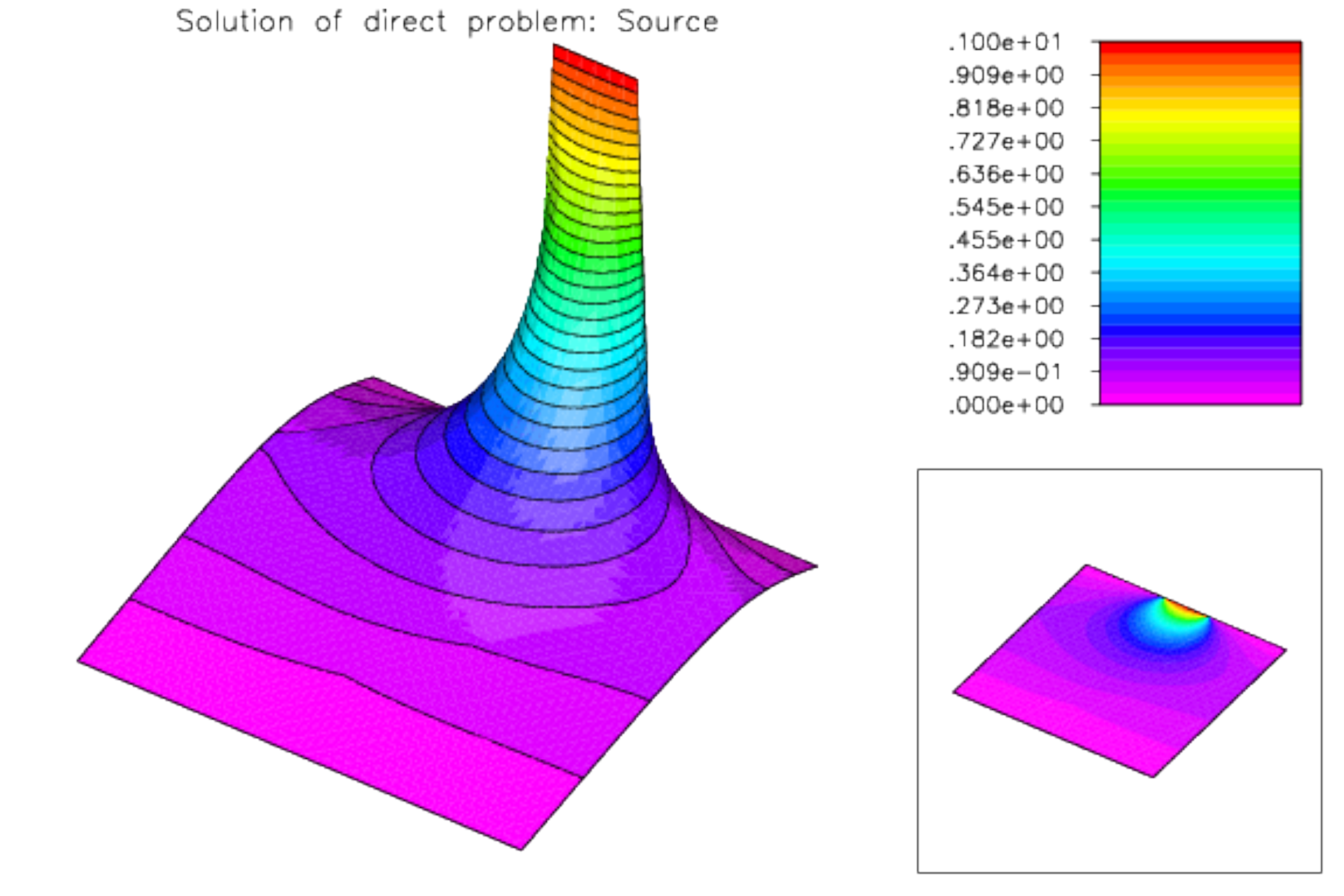} }
\centerline{\hfil (a) \hskip5cm  (b) \hfil}
\caption{\small Picture (a) show the doping profiles to be reconstructed in
the numerical experiments.} \label{fig:exsol-source}
\end{figure}

For this experiment concerning pointwise measurements of the current
density, we assume that only one measurement is available, i.e. $N = 1$
in (\ref{eq:data-bp-pm}). What concerns the numerical implementation,
we applied an iterative method of level set type to solve the identification
problem for $\gamma$ in (\ref{eq:num-d2nB}) (see \cite{LS03,FSL05,LMZ06}).
The domain $\Omega \subset \mathbb R^2$ is the unit square, and the boundary
parts are defined as follows
$$ \Gamma_1 \ := \  \{ (x,1) \, ;\ x \in (0,1) \} \, , \ \
   \Gamma_0 \ := \  \{ (x,0) \, ;\ x \in (0,1) \} \, , $$
$$ \partial\Omega_N \ := \ \{ (0,y) \, ;\ y \in (0,1) \} \cup 
   \{ (1,y) \, ;\ y \in (0,1) \} \, . $$
The fixed input $U$, is chosen to be a piecewise constant function
supported in $\Gamma_0$
$$  U(x) \ := \ \left\{ \begin{array}{rl}
      1, & |x - 0.5| \le h \\
      0, & {\rm else} \end{array} \right. $$
The doping profile to be reconstructed is shown in Figure~%
\ref{fig:exsol-source}~(a). In Figure~\ref{fig:exsol-source}~(b) the
voltage source $U$ (applied at $\Gamma_0$) and the corresponding solution
$\hat u$ of (\ref{eq:num-d2nB}) are shown.
In these pictures, as well as in the forthcoming ones, $\Gamma_1$ is the
lower left edge and $\Gamma_0$ is the top right edge (the origin corresponds
to the upper right corner).

\begin{figure}[t]
\centerline{ \includegraphics[width=10.0cm]{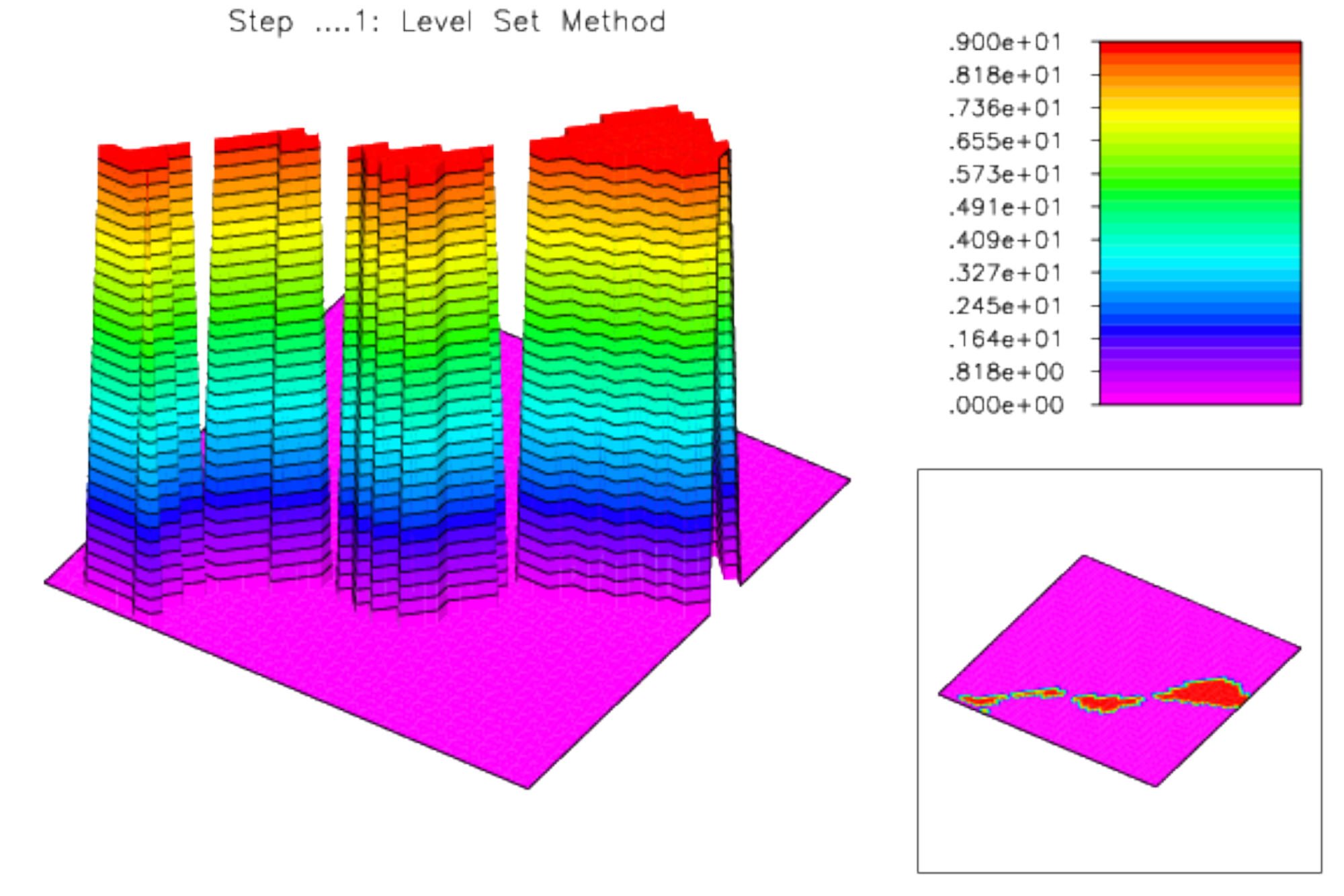} }
\centerline{ \includegraphics[width=10.0cm]{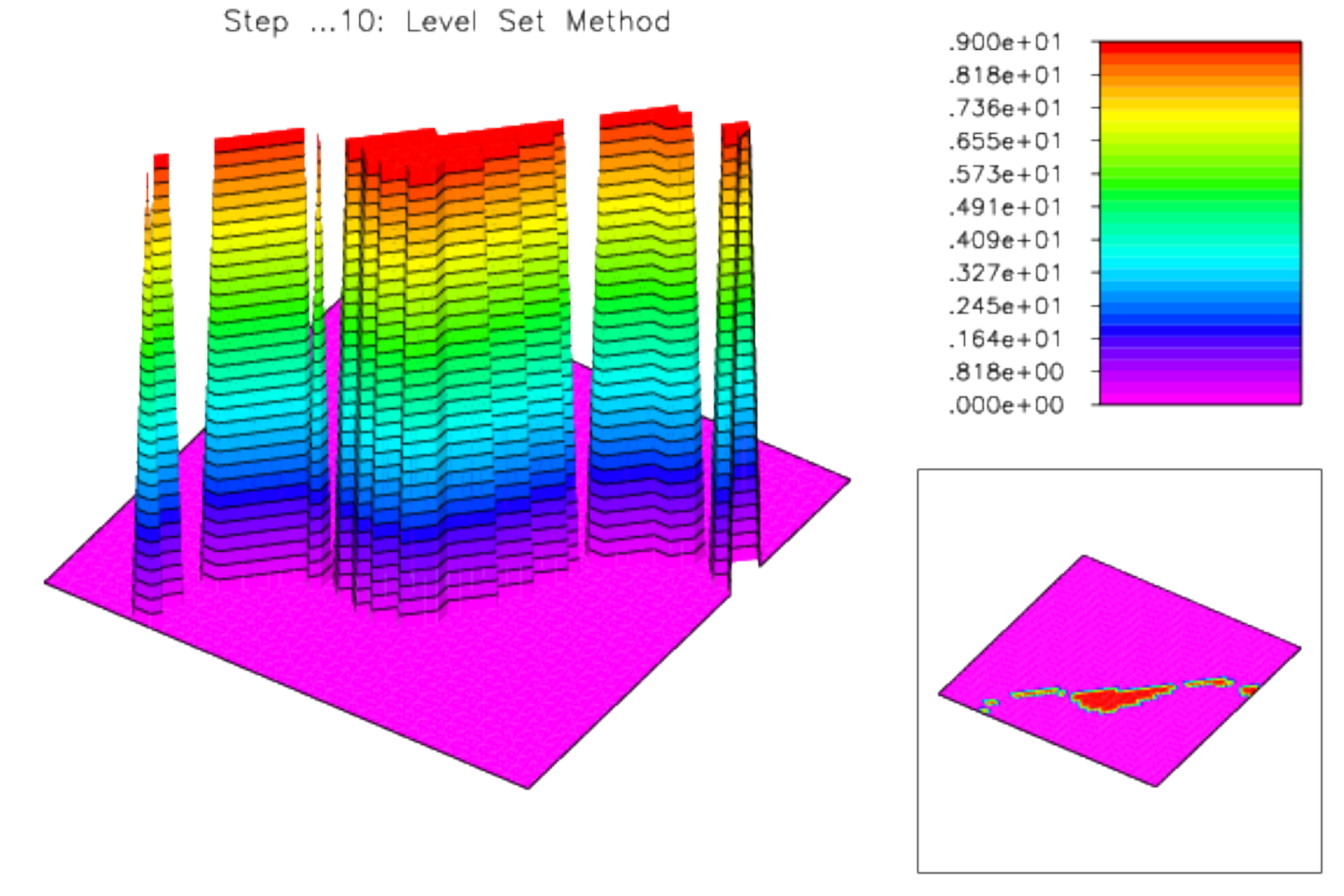} }
\centerline{ \includegraphics[width=10.0cm]{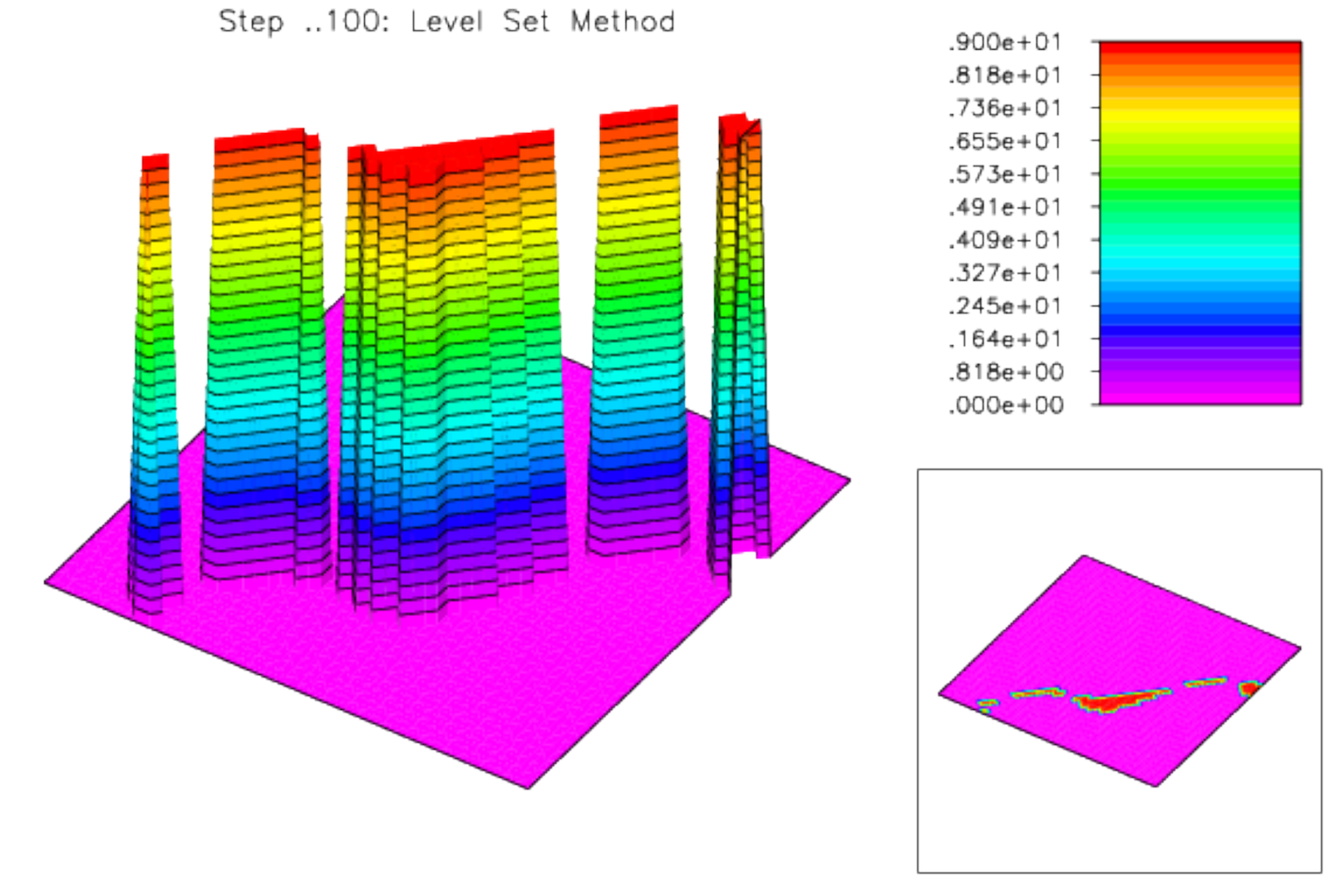} }
\caption{\small Experiment for the bipolar model with pointwise measurements
of the current density:
Reconstruction of the p-n junction in Figure~\ref{fig:exsol-source}~(a).
Evolution of the iteration error for exact data and one measurement of
the DN map $\Lambda_\gamma$ (i.e. $N = 1$ in (\ref{eq:data-bp-pm})).}
\label{fig:exp-bp-pm}
\end{figure}

In Figure~\ref{fig:exp-bp-pm} we present a numerical experiment for the
bipolar model with pointwise measurements of the current density. Here
exact data is used for the reconstruction of the p-n junction in Figure~%
\ref{fig:exsol-source}~(a). The pictures show plots of the iteration error
after 1, 10 and 100 steps of the level set method respectively.

What concerns the quality of the reconstruction of the P-N junction, the
level set approach considered in this paper brings much better results
than the Landweber-Kaczmarz approach implemented in \cite{BELM04}.
A possible explanation for the different performance of these methods is
the fact that the Landweber-Kaczmarz approach does not take into account
the assumption that the coefficient $\gamma$ in (\ref{eq:num-d2nB}) for
such application is a piecewise constant function.
The Landweber-Kaczmarz method tries to identify a real function defined on
$\Omega$, which is a much more complicated object than the original unknown
curve (the P-N junction).
Due to the nature of the level set approach, it incorporates in a natural
way the assumption that $\gamma$ is piecewise constant in $\Omega$.

\section*{Appendix}

Properties of silicon at room temperature

\begin{table}[ht]
\begin{center} \begin{tabular}{cl}
\hline \\[-2.3ex]
{ Parameter} & { \hfil Typical value \hfil} \\[0.3ex]
\hline \\[-2.0ex]
$\epsilon$   & $11.9 \ \epsilon_0$ \\
$\mu_n$      & $\approx 1500 \ {\rm cm}^2 \ {\rm V}^{-1} \ {\rm s}^{-1}$ \\
$\mu_p$      & $\approx  450 \ {\rm cm}^2 \ {\rm V}^{-1} \ {\rm s}^{-1}$ \\
$C_n$    & $2.8 \times 10^{-31} \ {\rm cm}^6 {\rm / s}$ \\
$C_p$    & $9.9 \times 10^{-32} \ {\rm cm}^6 {\rm / s}$ \\
$\tau_n$ & $10^{-6}\, {\rm s}$ \\
$\tau_p$ & $10^{-5}\, {\rm s}$ \\
\hline
\end{tabular} \end{center}
\caption{Typical values of main the constants in the model.
\label{tab:typ-val}}
\end{table}

\noindent
Relevant physical constants:

\begin{itemize}
\item[] Permittivity in vacuum:
       $\epsilon_0 = 8.85 \times 10^{-14} {\rm As \, V}^{-1} \, {\rm cm}^{-1}$;
\item[] Elementary charge:
       $q = 1.6 \times 10^{-19} {\rm As}$.
\end{itemize}

\subsection*{Acknowledgments}
The author was partially supported by the Brazilian National Research
Council CNPq, grants 305823/03-5 and 478099/04-5.


\end{document}